\magnification=1200
\documentstyle{amsppt}
\def\sign{\text{sign}\,}
\def\tr{\text{Tr}}
\def\prphi{\varphi_n(z)\varphi_{n+\ell}^{-1}(z)
       [\varphi_{n+\ell}^*(z)]^{-1} \varphi_n^*(z)}
\topmatter
\title
Rakhmanov's theorem for orthogonal matrix polynomials on the unit circle
\endtitle
\thanks
This research was supported by INTAS Research Network 03-51-6637, 
FWO Research Project G.0455.04 and K.U.Leuven research grant OT/04/21.
\endthanks
\author
Walter Van Assche
\endauthor
\affil
Katholieke Universiteit Leuven
\endaffil
\address
Department of Mathematics,
Katholieke Universiteit Leuven,
Celestijnenlaan 200 B,
B-3001 Leuven,
Belgium
\endaddress
\email
walter\@wis.kuleuven.be  \hfill\break
\indent {\it WWW home}:
http://www.wis.kuleuven.be/wis/analyse/walter.html
\endemail
\rightheadtext{Rakhmanov's theorem for matrix orthogonal polynomials}
\leftheadtext{Walter Van Assche}
\abstract
Rakhmanov's theorem for orthogonal polynomials on the unit circle gives a 
sufficient condition on the orthogonality measure for orthogonal
polynomials on the unit circle, in order that the reflection coefficients
(the recurrence coefficients in the Szeg\H{o} recurrence relation)
converge to zero. In this paper we give the analog for orthogonal
matrix polynomials on the unit circle.
\endabstract

\endtopmatter

\document

\head
1. Rakhmanov's theorem in the scalar case
\endhead

Let $\varphi_n(z) = \kappa_n z^n + \cdots$ $(n=0,1,2,\ldots)$,
with $\kappa_n > 0$,  be orthonormal polynomials on the
unit circle with respect to some positive measure $\mu$:
$$   \int_0^{2\pi} \varphi_n(z) \overline{\varphi_m(z)} \, d\mu(\theta)
        = \delta_{m,n},  \qquad z= e^{i\theta}. $$
We denote the monic polynomials by $\Phi_n(z) = \varphi_n(z)/\kappa_n$.
These monic polynomials satisfy a useful recurrence relation
$$  \Phi_n(z) = z \Phi_{n-1}(z) + \Phi_n(0) \Phi_{n-1}^*(z), \tag 1.1 $$
where $\Phi_n^*(z) = z^n \overline{\Phi}_n(1/z)$ is the reversed polynomial
(see, e.g., \cite{12,15}).
The coefficients $\Phi_n(0)$, which act as recurrence coefficients in this
recurrence relation, are known as reflection coefficients and 
$\alpha_n = -\overline{\Phi_{n+1}(0)}$ are called Verblunsky coefficients in \cite{12}. 
It is well known
that all the zeros $z_{k,n}$ of $\varphi_n$ lie in
the open unit disk, and hence $|\Phi_n(0)| = \prod_{k=1}^n |z_{k,n}| < 1$.
Moreover, 
$$   \frac{\kappa_{n-1}^2}{\kappa_n^2} = 1 - |\Phi_n(0)|^2,  \tag 1.2 $$
so that the reflection coefficients allow us to compute the monic orthogonal
polynomials recursively using \thetag{1.1}, but also the orthonormal
polynomials using in addition \thetag{1.2}. Conversely, given a sequence
$a_n$ of complex numbers for which $|a_n| < 1$ for all $n > 0$, then
polynomials satisfying
$$   \Phi_n(z) = z \Phi_{n-1}(z) + a_n \Phi_{n-1}^*(z), $$
with $\Phi_0 = 1$, will be monic orthogonal polynomials on the unit circle
for a unique orthogonality measure $\mu$, and the reflection
coefficients of these polynomials are $\Phi_n(0) = a_n$ (Favard's theorem
for orthogonal polynomials on the unit circle or Geronimus's theorem, see, e.g.,  \cite{3} and \cite{12}).

It is straightforward to see that the condition
$$      \lim_{n\to \infty} \Phi_n(0) = 0        \tag 1.3         $$
implies that
$$  \lim_{n \to \infty} \frac{\kappa_{n-1}^2}{\kappa_n^2} = 1  \tag 1.4 $$
and from \thetag{1.1} we also see that
$$  \lim_{n \to \infty} \frac{\Phi_n(z)}{\Phi_{n-1}(z)} = z \tag 1.5  $$
uniformly for $z$ on the unit circle ${\Bbb T} = \{ z \in {\Bbb C}: |z|=1\}$,
because on the unit circle we have $|\Phi_n^*(z)/\Phi_n(z)|=1$. Combined
with the ratio behavior \thetag{1.4} this also gives
ratio behavior of the orthonormal polynomials
$$  \lim_{n \to \infty} \frac{\varphi_n(z)}{\varphi_{n-1}(z)} = z ,
 \qquad z \in {\Bbb T}.    \tag 1.6  $$
This indicates that the orthonormal polynomials $\varphi_n(z)$
behaves very much like the polynomials $z^n$, which are the
orthonormal polynomials on the unit circle for Lebesgue measure
$d\theta/2\pi$.
The condition \thetag{1.3} is therefore a very natural
condition for various asymptotic properties of the orthogonal
polynomials. For this reason it is of interest to find conditions on the
orthogonality measure $\mu$ that imply \thetag{1.3}. Rakhmanov \cite{11}
(see also \cite{8})
proved that a rather mild condition on the size of $\mu$ on the unit
circle is sufficient.

\proclaim{Rakhmanov's Theorem}
Suppose that $\mu' > 0$ almost everywhere on ${\Bbb T}$. Then
$\lim_{n \to \infty} \Phi_n(0) = 0$.
\endproclaim

This condition is not necessary: we know examples of discrete measures
and singularly continuous measures, hence
measure without absolutely continuous component, for which \thetag{1.3} holds \cite{6}, 
\cite{7}, \cite{16}, \cite{14}. 
A nice proof of Rakhmanov's theorem uses the following two equivalences
(Nevai \cite{9} \cite{10}, Li and Saff \cite{5})
$$   \lim_{n \to \infty} \Phi_n(0) = 0 
     \Longleftrightarrow
    \lim_{n \to \infty} \inf_{\ell \geq 1} \frac{1}{2\pi}
   \int_0^{2\pi} \left| \frac{|\varphi_n(z)|^2}{|\varphi_{n+\ell}(z)|^2}
           - 1 \right| \, d\theta = 0 ,   
      \tag 1.7         $$
and
$$   \mu'(\theta) > 0 \text{ a.e. on } [0,2\pi) 
   \Longleftrightarrow
  \lim_{n \to \infty} \sup_{\ell \geq 1} \frac{1}{2\pi}
   \int_0^{2\pi} \left| \frac{|\varphi_n(z)|^2}{|\varphi_{n+\ell}(z)|^2}
           - 1 \right| \, d\theta = 0 .  
         \tag 1.8        $$

In this paper we will investigate the analog of Rakhmanov's theorem
for orthogonal matrix polynomials on the unit circle. In Section 2
we will introduce the necessary background on orthogonal matrix
polynomials and present the matrix analogs of the Szeg\H{o}
recurrence \thetag{1.1} and the reflection coefficients, which turn
out to be matrices. In Section 3 we prove the matrix analog of the
characterization \thetag{1.7} and in Section 4 we deal with the
characterization \thetag{1.8}.  
 
\head
2. Orthogonal matrix polynomials on the unit circle
\endhead
 
Let $\rho$ be a $p\times p$ Hermitian matrix-valued measure on the unit
circle ${\Bbb T}$, then $\rho$ induces two matrix inner products. The
left inner product
$$  \langle P,Q \rangle_L = \int_0^{2\pi} P(z)\, d\rho(\theta)\, Q(z)^*,
\qquad z=e^{i\theta}, $$
where from now on $Q^*$ is the Hermitian transpose of the matrix $Q$, i.e.,
$Q^* = \overline{Q^T}$, and the right inner product
$$   \langle P,Q \rangle_R = \int_0^{2\pi} P(z)^*\, d\rho(\theta)\, Q(z),
\qquad z=e^{i\theta}. $$
Both inner products give rise to a sequence of orthonormal polynomials,
which we call left orthonormal polynomials $\varphi_n^L$ and right orthonormal polynomials $\varphi_n^R$
respectively:
$$  \int_0^{2\pi} \varphi_n^L(z)\, d\rho(\theta)\, [\varphi_m^L(z)]^*
    = \delta_{m,n}
  = \int_0^{2\pi} [\varphi_n^R(z)]^*\, d\rho(\theta)\, \varphi_m^R(z).
$$
These orthonormal polynomials are determined up to multiplication
by a unitary matrix (on the left for the left orthonormal polynomials
and on the right for the right orthonormal polynomials). For background and
general results on these orthogonal polynomials on the unit circle, we refer to \cite{1} and \cite{13}.

We use the notation $\tilde{P}(z) = z^n P(1/\bar{z})^*$ for the reversed matrix polynomial, whenever $P$ is a polynomial of degree at most $n$.
If we write the orthonormal polynomials as
$$  \align
  \varphi_n^L(z) &= L_{n,n} z^n + L_{n,n-1} z^{n-1} + \cdots + L_{n,0} \\
  \varphi_n^R(z) &= K_{n,n} z^n + K_{n,n-1} z^{n-1} + \cdots + K_{n,0}
  \endalign 
$$
then the leading coefficients $K_{n,n}$ and $L_{n,n}$ are non-singular
$p\times p$ matrices.
The identity $\langle \varphi_n^R,\tilde{\varphi}_n^L\rangle_R
= \langle \tilde{\varphi}_n^R,\varphi_n^L \rangle_L$ implies
$$   L_{n,0}^* L_{n,n} = K_{n,n} K_{n,0}^*, $$
and we can then introduce the reflection coefficients as 
$$  H_n = (L_{n,n}^*)^{-1} K_{n,0} = L_{n,0} (K_{n,n}^*)^{-1}. \tag 2.1 $$
One can show that $I-H_n H_n^*$ and $I - H_n^*H_n$ are positive definite,
so that $\|H_n\|_2 < 1$ for all $n > 0$, 
and
$$  \align
   (I-H_n^*H_n)^{1/2} &= K_{n,n}^{-1}K_{n-1,n-1}   \\
   (I-H_nH_n^*)^{1/2} &= (L_{n,n}^*)^{-1}L_{n-1,n-1}^*,
  \endalign $$
which is the matrix analog of \thetag{1.2}, and furthermore we have the
recurrences
$$ \align
   (I-H_n H_n^*)^{1/2} \varphi_n^L(z) &= z \varphi_{n-1}^L(z)
     + H_n \tilde{\varphi}_{n-1}^R(z)  \tag 2.2 \\
    \varphi_n^R(z) (I-H_n^*H_n)^{1/2} &= z \varphi_{n-1}^R(z)
    + \tilde{\varphi}_{n-1}^L(z)H_n ,  \tag 2.3 
  \endalign $$
which are the matrix analogs of the Szeg\H{o} recurrences \thetag{1.1},
but for orthonormal polynomials. Observe that
for orthogonal matrix polynomials on the unit circle we need to work with
both left and right orthogonal polynomials. Again there is a converse
result (Favard's theorem for orthogonal matrix polynomials on the
unit circle) saying that matrix polynomials $\varphi_n^L$ and
$\varphi_n^R$ satisfying the recurrences \thetag{2.2}--\thetag{2.3},
with initial conditions $\varphi_0^L$ and $\varphi_0^R$ which are
non-singular matrices with 
$[\varphi_0^L]^*\varphi_0^L = \varphi_0^R [\varphi_0^R]^*$ and
with matrix coefficients $H_n$ for which $\|H_n\|_2 < 1$ for all $n > 0$,
are always orthonormal matrix polynomials on the unit circle for some
positive matrix measure $\rho$ \cite{1, Thm.~15 on p.~155}.

Using \thetag{2.2} and \thetag{2.3} and the identity 
$(I-AA^*)^{-1}A = A(1-A^*A)^{-1}$, we can find
$$  \tilde{\varphi}_n^L(z)\varphi_n^L(\xi) - \xi \tilde{\varphi}_{n-1}^L(z)
   \varphi_{n-1}^L(\xi) = \varphi_n^R(z)\tilde{\varphi}_n^R(z)
   - z \varphi_{n-1}^R(z)\tilde{\varphi}_{n-1}^R(\xi). $$
Multiply both sides by $z^{-n}$, then summing over $n$ gives
$$  z^{-n} \varphi_n^R(z) \tilde{\varphi}_n^R(\xi)
   = (1-\xi/z) \sum_{k=0}^n z^{-k}\tilde{\varphi}_k^L(z) \varphi_k^L(\xi)
     + \xi z^{-n} \tilde{\varphi}_n^L(z) \varphi_n^L(\xi). $$
Now use $\tilde{\varphi}_n^L(z) = z^n \varphi_n^L(1/\bar{z})^*$ and
replace $z$ by $1/\bar{z}$, then this gives the {\it Christoffel-Darboux 
formula}
$$  (1-\xi \bar{z}) \sum_{k=0}^n [\varphi_k^L(z)]^* \varphi_k^L(\xi)
     = [\tilde{\varphi}_n^R(z)]^* \tilde{\varphi}_n^R(\xi) - \xi \bar{z}
   [\varphi_n^L(z)]^*  \varphi_n^L(\xi). \tag 2.4   $$
In a similar way we can obtain the dual formula
$$  (1-\xi \bar{z}) \sum_{k=0}^n \varphi_k^R(\xi) [\varphi_k^R(z)]^*
     = \tilde{\varphi}_n^L(\xi) [\tilde{\varphi}_n^L(z)]^*
        - \xi \bar{z} \varphi_n^R(\xi) [\varphi_n^R(z)]^*. \tag 2.5 $$
An important consequence of this Christoffel-Darboux formula
is when we take $z = \xi$ on the unit circle, which gives the identity
$$  \varphi_n^R(z) \tilde{\varphi}_n^R(z) =
   \tilde{\varphi}_n^L(z) \varphi_n^L(z),\qquad z=e^{i\theta}. \tag 2.6   $$
Observe that on the unit circle $\tilde{P}(z) = z^n P(z)^*$ so that
$z^{-n}\tilde{\varphi}_n^L(z) \varphi_n^L(z)$ is a positive
definite matrix for every $z$ on the unit circle. Of particular
use is the absolutely continuous matrix measure
$$   d\rho_n(\theta)
  = \left( [\varphi_n^L(z)]^* \varphi_n^L(z) \right)^{-1}\, d\theta/2\pi
  = \left( \varphi_n^R(z) [\varphi_n^R(z)]^* \right)^{-1}\, d\theta/2\pi,
   \qquad z = e^{i\theta},  \tag 2.7 $$
because the first $n+1$ orthogonal matrix polynomials are also orthogonal
with respect to this measure (\cite{1, Eq. (71) on p.~155})
$$  \frac{1}{2\pi} \int_0^{2\pi} \varphi_k^L(z)\ 
    \left( [\varphi_n^L(z)]^* \varphi_n^L(z) \right)^{-1}\
    [\varphi_m^L(z)]^*\, d\theta = \delta_{k,m}, \qquad k,m \leq n,
                                  \tag 2.8  $$
and similarly for the right orthogonal matrix polynomials.

We are now ready to state Rakhmanov's theorem for orthogonal matrix
polynomials on the unit circle.

\proclaim{Theorem}
Suppose that $\rho$ is a matrix measure on the unit circle with
$d\rho(\theta) = \rho'(\theta) \, d\theta/2\pi + d\rho_s(\theta)$
and $\rho_s$ the singular part of the measure. If $\det \rho'(\theta) > 0$
almost everywhere on $[0,2\pi)$ then $\lim_{n \to \infty} H_n = 0$.
\endproclaim

We will prove this using appropriate matrix analogs of \thetag{1.7}
(see section 3) and \thetag{1.8} (see section 4).

\head
3. Characterization of $\lim_{n \to \infty} H_n = 0$
\endhead

\proclaim{Lemma 1}
Let $\varphi_n^L(z)$ and $\varphi_n^R(z)$ be left and right orthonormal
matrix polynomials with reflection coefficients $H_n$, then
$$  \multline \lim_{n \to \infty} H_n = 0 
         \Longleftrightarrow \\
    \lim_{n \to \infty} \inf_{\ell \geq 1}\frac{1}{2\pi} \int_0^{2\pi} 
  \| \varphi_n^L(z) \varphi_{n+\ell}^L(z)^{-1}
    [\varphi_{n+\ell}^L(z)^*]^{-1} \varphi_n^L(z)^* - I \|_2 \,
    d\theta = 0 .  
   \endmultline  \tag 3.1 $$
\endproclaim

\demo{Proof of $\Rightarrow$}
The recurrence relation \thetag{2.2} gives
$$   (I-H_n H_n^*)^{1/2} \varphi_n^L(z) \varphi_{n-1}^L(z)^{-1}
    = z I + H_n \tilde{\varphi}_{n-1}^R(z) \varphi_{n-1}^L(z)^{-1},
   \tag 3.2 $$
and similarly \thetag{2.3} gives
$$    \varphi_{n-1}^R(z)^{-1} \varphi_n^R(z)(I-H_n^*H_n)^{1/2}
    = z I +  \varphi_{n-1}^R(z)^{-1} \tilde{\varphi}_{n-1}^L(z) H_n.
  \tag 3.3 $$
Observe that $\tilde{\varphi}_{n-1}^R(z) \varphi_{n-1}^L(z)^{-1} =
\varphi_{n-1}^R(z)^{-1} \tilde{\varphi}_{n-1}^L(z)$ whenever $z$ is on
the unit circle, which follows from \thetag{2.6}. Furthermore
this rational matrix function is a unitary matrix when $z$ is on the
unit circle. To see this, we examine the product
$$   [\tilde{\varphi}_{n-1}^R(z) \varphi_{n-1}^L(z)^{-1}]
   [\tilde{\varphi}_{n-1}^R(z) \varphi_{n-1}^L(z)^{-1}]^*, $$
which is equal to 
$$  \tilde{\varphi}_{n-1}^R(z) [\varphi_{n-1}^L(z)^*
   \varphi_{n-1}^L(z)]^{-1} \tilde{\varphi}_{n-1}^R(z)^*. $$
Use \thetag{2.6} to replace the product of the two matrix polynomials, then
we find that this expression is 
$$  \tilde{\varphi}_{n-1}^R(z) [\varphi_{n-1}^R(z)
   \varphi_{n-1}^R(z)^*]^{-1} \tilde{\varphi}_{n-1}^R(z)^* $$
and this is
$$   \tilde{\varphi}_{n-1}^R(z) [\varphi_{n-1}^R(z)^*]^{-1}\ 
   \varphi_{n-1}^R(z)^{-1} \tilde{\varphi}_{n-1}^R(z)^*. $$
On the unit circle we have $\tilde{\varphi}_{n-1}^R(z) = 
z^{n-1} \varphi_{n-1}^R(z)^*$, hence the expression reduces to the unit matrix, so that
$$   [\tilde{\varphi}_{n-1}^R(z) \varphi_{n-1}^L(z)^{-1}]
   [\tilde{\varphi}_{n-1}^R(z) \varphi_{n-1}^L(z)^{-1}]^* = I $$
and $\tilde{\varphi}_{n-1}^R(z) \varphi_{n-1}^L(z)^{-1}$ is unitary.
Consequently for the spectral norm we have
$\|\tilde{\varphi}_{n-1}^R(z) \varphi_{n-1}^L(z)^{-1}\|_2 = 1$.
Returning to \thetag{3.2}, we now have
$$  \|(I-H_n H_n^*)^{1/2} \varphi_n^L(z) \varphi_{n-1}^L(z)^{-1}
      - z I \|_2 \leq \|H_n\|_2 $$
and similarly \thetag{3.3} gives
$$  \|\varphi_{n-1}^R(z)^{-1} \varphi_n^R(z)(I-H_n*H_n)^{1/2} - zI \|_2
    \leq  \|H_n\|_2. $$
Hence $\lim_{n \to \infty} H_n = 0$ implies
$$  \lim_{n \to \infty} \varphi_n^L(z) \varphi_{n-1}^L(z)^{-1} = z I ,
      \quad
  \lim_{n \to \infty} \varphi_{n-1}^R(z)^{-1} \varphi_n^R(z)= zI ,   $$
uniformly for $z \in {\Bbb T}$. This implies that for $\ell \geq 1$ fixed
$$ \lim_{n \to \infty}  \varphi_{n}^L(z) \varphi_{n+\ell}^L(z)^{-1} =
    z^{-\ell} I,   $$
uniformly on the unit circle, which implies that
$$  \lim_{n \to \infty} \inf_{\ell \geq 1}\frac{1}{2\pi} \int_0^{2\pi} 
  \| \varphi_n^L(z) \varphi_{n+\ell}^L(z)^{-1}
     [\varphi_{n+\ell}^L(z)^*]^{-1} \varphi_n^L(z)^* - I \|_2 \,
    d\theta = 0,  $$
which is what we wanted to prove. Observe that we also get a similar
result for the right orthogonal polynomials.

\enddemo

\demo{Proof of $\Leftarrow$}
Here we rely on the identity
$$   \langle z \varphi_n^L,  \tilde{\varphi}_n^R \rangle_L = -H_{n+1}.  
    \tag 3.4   $$
Indeed, if we use \thetag{2.2} (replacing $n$ by $n+1$ everywhere), then
$$  \align
   \langle z\varphi_n^L, \tilde{\varphi}_n^R \rangle_L		
   &= \int_0^{2\pi} z\varphi_n^L(z)\, d\rho(\theta)\,        [\tilde{\varphi}_n^R(z)]^* \\
   &= (I - H_{n+1}H_{n+1}^*)^{1/2} \langle     \varphi_{n+1}^L,\tilde{\varphi}_n^R \rangle_L
  - H_{n+1} \langle \tilde{\varphi}_n^R , \tilde{\varphi}_n^R \rangle_L \\
   &=  - H_{n+1} \langle \tilde{\varphi}_n^R , \tilde{\varphi}_n^R \rangle_L
 \endalign $$
where the last step follows from the orthogonality. Now
$\langle \tilde{P_n},\tilde{Q_n} \rangle_L = \langle P_n,Q_n \rangle_R$,
hence
$$   \langle z\varphi_n^L, \tilde{\varphi}_n^R \rangle_L
    = -H_{n+1} \langle \varphi_n^R, \varphi_n^R \rangle_R = -H_{n+1} $$
where we used the orthonormality. This shows that \thetag{3.4}
indeed holds.

 From the finite orthogonality \thetag{2.8} we can easily deduce the
following result for the measures $\rho$ and $\rho_n$, given by \thetag{2.7}:
$$  \int_0^{2\pi} P_k(z)\, d\rho(\theta)\, [Q_m(z)]^* =
     \frac{1}{2\pi} \int_0^{2\pi} P_k(z)\, 
  ([\varphi_n^L(z)]^* \varphi_n^L(z))^{-1}\, Q_m(z)\, d\theta,  $$
for all matrix polynomials $P_k$ and $Q_m$ of degree $k\leq n$ and $m\leq n$
respectively, by expanding $P_k$ and $Q_m$ in a Fourier series using
the left orthonormal polynomials. Since $z \varphi_n^L(z)$ is a matrix
polynomial of degree $n+1$ and $\tilde{\varphi}_n^R(z)$ is of degree $n$,
we therefore have from \thetag{3.4}
$$  -H_{n+1} = \frac{1}{2\pi} \int_0^{2\pi}
    z \varphi_n^L(z)\,  ([\varphi_{n+\ell}^L(z)]^* 
     \varphi_{n+\ell}^L(z))^{-1}\, [\tilde{\varphi}_n^R(z)]^*\, d\theta, $$
for every $\ell \geq 1$.
We can write this integral as
$$ \multline
 \frac{1}{2\pi} \int_0^{2\pi}
    \varphi_n^L(z) \varphi_{n+\ell}^L(z)^{-1} \ 
     [\varphi_{n+\ell}^L(z))^*]^{-1} \varphi_n^L(z)^*\
  z [\varphi_n^L(z)^*]^{-1} \tilde{\varphi}_n^R(z)]^*\, d\theta \\
   = \frac{1}{2\pi} \int_0^{2\pi}
     [\varphi_n^L(z)\varphi_{n+\ell}^L(z)^{-1}]\, 
     [\varphi_n^L(z)\varphi_{n+\ell}^L(z)^{-1}]^*  \     z   \tilde{\varphi}_n^L(z)^{-1} \varphi_n^R(z)\, d\theta. 
 \endmultline $$
By the calculus of residues we also have
$$ \frac{1}{2\pi} \int_0^{2\pi}
    z \tilde{\varphi}_n^L(z)^{-1} \varphi_n^R(z)\,  d\theta = 0, $$
hence we get the formula
$$ -H_{n+1} =  \frac{1}{2\pi} \int_0^{2\pi}
     \left( [\varphi_n^L(z)\varphi_{n+\ell}^L(z)^{-1}]\, 
     [\varphi_n^L(z)\varphi_{n+\ell}^L(z)^{-1}]^* - I \right)\,  z        \tilde{\varphi}_n^L(z)^{-1} \varphi_n^R(z)\, d\theta .
  $$
Recall that $\tilde{\varphi}_n^L(z)^{-1} \varphi_n^R(z)$ is unitary
whenever $z \in {\Bbb T}$, so that by taking the spectral norm we get
$$ \|H_{n+1}\|_2 \leq \frac{1}{2\pi} \int_0^{2\pi}
     \| [\varphi_n^L(z)\varphi_{n+\ell}^L(z)^{-1}]\, 
     [\varphi_n^L(z)\varphi_{n+\ell}^L(z)^{-1}]^* - I \|_2 \, d\theta. $$
Hence if 
$$  \lim_{n \to \infty} \inf_{\ell \geq 1} \frac{1}{2\pi} \int_0^{2\pi}
     \| [\varphi_n^L(z)\varphi_{n+\ell}^L(z)^{-1}]\, 
   [\varphi_n^L(z)\varphi_{n+\ell}^L(z)^{-1}]^* - I \|_2 \, d\theta = 0, $$
the obviously $H_n \to 0$, which is what we wanted to prove. \qed    
\enddemo

\head
4. Characterization of $\det \rho'(\theta) > 0$ almost everywhere
\endhead

In this section we will only use left orthogonal polynomials $\varphi_n^L$
and to simplify the notation we therefore will drop the superscript $L$.

In general the matrix measure $\rho$ will consist of an absolutely
continuous part with Radon-Nikodym derivative $\rho'$, and a
singular part $\rho_s$. A remarkable fact is that when 
$\det \rho'(\theta) > 0$ almost everywhere, then the singular part
$\rho_s$ does not interfere in the ratio asymptotic behavior.
We can indeed annihilate the singular part using the same ideas
as in  \cite{10}.

\proclaim{Lemma 2}
Suppose $\rho_s$ is a positive definite matrix measure on the unit circle
which is singular with respect to the Lebesgue matrix measure.
Then there exists a sequence of matrix functions $G_n$ on the unit circle, such that
$G_n G_n^* \leq I$,
$$   \lim_{n \to \infty} G_n G_n^* = I, \qquad \text{almost everywhere
on } {\Bbb T},  $$ 
and 
$$  \lim_{n \to \infty} \int_0^{2\pi} G_n(\theta) \, d\rho_s(\theta)\,
     G_n^*(\theta) = 0. $$
\endproclaim

\demo{Proof}
 From Lemma 5 in \cite{10} we know that if $\mu_s$ is a (scalar) singular measure
on the unit circle, then there is a sequence of real-valued
$2\pi$-periodic continuous functions
$h_n$ on the unit circle, such that $0 \leq h_n(\theta) \leq 1$
for all $\theta \in [0,2\pi)$, with
$$  \lim_{n \to \infty} h_n(\theta) = 1, \qquad \text{almost everywhere} $$
and
$$  \lim_{n \to \infty} \int_0^{2\pi} h_n(\theta)\, d\mu_s(\theta) = 0. $$
Let
$$  \rho_s(\theta) = V(\theta) \pmatrix
			\rho_{s,1}(\theta) &  &  & 0 \\
			  & \rho_{s,2}(\theta) & &  \\
			  &  & \ddots & \\
			  0&  &  & \rho_{s,p}(\theta)
			\endpmatrix  V^*(\theta)  $$
be the Schur decomposition of $\rho_s$, with $V$ a unitary matrix 
function. Each eigenvalue $\rho_{s,i}$ is singular with respect to Lebesgue measure, hence there exists a sequence $h_{n,i}$ of 
$2\pi$-periodic continuous functions such that 
$0 \leq h_{n,i}(\theta) \leq 1$, with
$$  \lim_{n \to \infty} h_{n,i}(\theta) = 1 \qquad \text{almost everywhere} $$
and
$$  \lim_{n \to \infty} \int_0^{2\pi} h_{n,i}(\theta) \, d\rho_{s,i}(\theta) = 0. $$
Consider the sequence of matrix functions
$$   G_n = \pmatrix
         h_{n,1}^{1/2}  & & & 0 \\
	    & h_{n,2}^{1/2} & & \\			
          & & \ddots & \\
         0 & & & h_{n,p}^{1/2} 
		\endpmatrix   V^*(\theta), $$
then $G_nG_n^* \leq I$, 
$$  \lim_{n \to \infty} G_nG_n^* = I, \qquad \text{almost everywhere} $$
and
$$  \lim_{n \to \infty} \int_0^{2\pi} G_n \, d\rho_s\, G_n^* = 0, $$
which is what we wanted to proof. \qed
\enddemo

\proclaim{Lemma 3}
Let $\varphi_n$ be the left  orthonormal
matrix polynomials with reflection coefficients $H_n$, then
$$  \multline  \det \rho'(x) > 0 \text{ almost everywhere} 
   \Longleftrightarrow \\
    \lim_{n \to \infty} \sup_{\ell \geq 1}\frac{1}{2\pi} \int_0^{2\pi} 
  \| \varphi_n (z) \varphi_{n+\ell} (z)^{-1}
  [\varphi_{n+\ell} (z)^*]^{-1} \varphi_n (z)^* - I \|_2 \,
    d\theta = 0 .  
  \endmultline     $$
\endproclaim

\demo{Proof of $\Leftarrow$}
First of all we observe that
$$  \frac{1}{2\pi} \int_0^{2\pi} P_m(z)\ \left( \varphi_n (z)^*
   \varphi_n (z) \right)^{-1}\ Q_m(z)^* \, d\theta
   = \int_0^{2\pi} P_m(\theta)\, d\rho(\theta)\, Q_m(\theta)^* $$
holds for all matrix polynomials $P_m$ and $Q_m$
of degree $m \leq n$. This follows from \thetag{2.8}. This means that
$$ \multline   \int_0^{2\pi} F_m(\theta) \varphi_n (z)\, \rho'(\theta)\,
       \varphi_n (z)^* G_m(\theta)^*\, d\theta 
   + \int_0^{2\pi} F_m(\theta) \varphi_n (z)\, d\rho_s(\theta)\,
       \varphi_n (z)^* G_m(\theta)^*  \\
   = \frac{1}{2\pi} \int_0^{2\pi} F_m(\theta) \varphi_n (z)
   \left( \varphi_{n+\ell} (z)^* \varphi_{n+\ell} (z) \right)^{-1}
     \varphi_n (z)^* G_m(\theta)^*\, d\theta, \qquad \ell \geq 2m , 
 \endmultline $$
for all trigonometric matrix polynomials $F_m$ and $G_m$ of degree
at most $m$. This gives 
$$  \multline \frac{1}{2\pi} \int_0^{2\pi}
    F_m(\theta) \left[ \varphi_n (z)\,     2\pi    \rho'(\theta)\,
       \varphi_n (z)^* - I \right] G_m(\theta)^*\, d\theta \\
   =  -  \int_0^{2\pi} F_m(\theta) \varphi_n (z)\, d\rho_s(\theta)\,
       \varphi_n (z)^* G_m(\theta)^* \\ 
 + \frac{1}{2\pi} \int_0^{2\pi} F_m(\theta)\left[ \varphi_n (z)
   \left( \varphi_{n+\ell} (z)^* \varphi_{n+\ell} (z) \right)^{-1}
     \varphi_n (z)^* - I \right] G_m(\theta)^*\, d\theta, 
  \endmultline $$
whenever $\ell \geq 2m$, and if we take the spectral norm, then
$$  \multline
\| \frac{1}{2\pi} \int_0^{2\pi} F_m(\theta)
 \left[ \varphi_n (z)\, 2\pi \rho'(\theta)\,
       \varphi_n (z)^* - I \right] G_m(\theta)^*\, d\theta\ \|_2 \\
   \leq \|\int_0^{2\pi} F_m(\theta) \varphi_n (z)\, d\rho_s(\theta)\,
       \varphi_n (z)^* G_m(\theta)^* \|_2 \\
+ \max_{\theta \in [0,2\pi)} \|F_m(\theta)\|_2
         \max_{\theta \in [0,2\pi)} \|G_m(\theta)\|_2 \\ \times
 \sup_{\ell \geq 1} \frac{1}{2\pi} \int_0^{2\pi} \| \varphi_n (z)
   \left( \varphi_{n+\ell} (z)^* \varphi_{n+\ell} (z) \right)^{-1}
     \varphi_n (z)^* - I \|_2   \, d\theta, 
\endmultline   $$
for all trigonometric matrix polynomials $P_m$ and $Q_m$. Every $2\pi$-periodic continuous matrix function can be uniformly approximated by
trigonometric matrix polynomials. Therefore we also have
$$ \multline
 \| \frac{1}{2\pi} \int_0^{2\pi} F(\theta) \left[ \varphi_n (z)\,  
  2\pi \rho'(\theta)\,
             \varphi_n (z)^* -I \right] G(\theta)^*\, d\theta\ \|_2 \\
   \leq \| \int_0^{2\pi} F(\theta) \varphi_n (z)\, d\rho_s(\theta)\,
       \varphi_n (z)^* G(\theta)^* \|_2 \\
 +  \max_{\theta \in [0,2\pi)} \|F(\theta)\|_2
         \max_{\theta \in [0,2\pi)} \|G(\theta)\|_2 \\ \times
 \sup_{\ell \geq 1} \frac{1}{2\pi} \int_0^{2\pi} \| \varphi_n (z)
   \left( \varphi_{n+\ell} (z)^* \varphi_{n+\ell} (z) \right)^{-1}
     \varphi_n (z)^* -I \|_2   \, d\theta,   
  \endmultline  \tag 3.8  $$
for all $2\pi$-periodic continuous matrix functions $F$ and $G$. Furthermore, every matrix function $F$ for which
$\sup_{\theta \in [0,2\pi)} \| F(\theta)\|_2 < \infty$ 
(i.e., $F \in L^\infty[0,2\pi]$) can be approximated
pointwise by $2\pi$-periodic continuous matrix functions $F_k$ with
$\sup_{\theta \in [0,2\pi)} \|F_k(\theta)\|_2 = 
\sup_{\theta \in [0,2\pi)} \|F(\theta)\|_2$, hence \thetag{3.8}
also holds for matrix functions $F, G \in L^\infty[0,2\pi]$.
A useful choice is to take $F=HP$ and $G=P$, where $P$ is a unitary
matrix function such that $P(\theta) \varphi_n (z) 2\pi \rho'(\theta)
\varphi_n (z)^* P(\theta)^*$ is a diagonal matrix $D_n$ containing
the eigenvalues $d_{1,n}(\theta), \ldots, d_{p,n}(\theta)$ of
the positive definite matrix $\varphi_n (z) 2\pi \rho'(\theta) 
\varphi_n (z)^*$, and $H$ is the unitary diagonal matrix with
entries $\sign (d_{1,n}-1), \sign (d_{2,n}-1), \ldots, \sign (d_{p,n}-1)$.
Then we get
$$ \multline  \max_{1 \leq k \leq p} \int_0^{2\pi} | d_{k,n}(\theta)-1|\,    d\theta \\
   \leq \sup_{\ell \geq 1}  \int_0^{2\pi} \| \varphi_n (z)
   \left( \varphi_{n+\ell} (z)^* \varphi_{n+\ell} (z) \right)^{-1}
     \varphi_n (z)^* - I\|_2   \, d\theta.  
\endmultline  \tag 3.9
$$
Observe that $\det \varphi_n (z) \rho'(\theta) \varphi_n (z)^*
= [\det \varphi_n (z)]^2 \det \rho'(\theta)$, and since $\varphi_n $
has no zeros on the unit circle, it follows that $\rho'(\theta) = 0$
if and only $d_{k,n}(\theta)=0$ for some $k=1,2,\ldots,p$. Hence
$$  A := \{ \theta \in [0,2\pi): \rho'(\theta) = 0 \} =
  \bigcup_{k=1}^p \{ \theta \in [0,2\pi): d_{k,n}(\theta)=0 \} :=
  \bigcup_{k=1}^p A_k . $$
Let $m$ be  Lebesgue measure on $[0,2\pi)$, i.e.,
$dm(\theta) = d\theta$, then
$$  m(A) =  \int_A 1\, d\theta 
        = \int_{\cup A_k} 1 \, d\theta 
        \leq \sum_{k=1}^p  \int_{A_k} 1\, d\theta. $$
Now, on $A_k$ we have $d_{k,n} = 0$, hence
$$  m(A) \leq \sum_{k=1}^p  \int_{A_k} |d_{k,n}(\theta)-1|\,
   d\theta 
   \leq \sum_{k=1}^p  \int_0^{2\pi} |d_{k,n}(\theta)-1|\,
   d\theta. $$
Now use \thetag{3.9} to conclude
$$  m(A) \leq p\ \sup_{\ell \geq 1}  \int_0^{2\pi} \| \varphi_n (z)
   \left( \varphi_{n+\ell} (z)^* \varphi_{n+\ell} (z) \right)^{-1}
     \varphi_n (z)^* - I\|_2   \, d\theta. $$
Hence if the right hand side tends to zero, then $\det \rho' > 0$ almost
everywhere on the unit circle. \qed
\enddemo

In order to prove the necessary part of the lemma, we will 
need some inequalities for the trace of matrices.
We will always be using $p\times p$ matrices and for
a matrix $A$ we will denote its eigenvalues by
$\lambda_i(A)$ $(i=1,2,\dots,p)$ and its singular values
by $\sigma_i(A)$ $(i=1,2,\ldots,p)$. Recall that for a
matrix $A$ we always have $\sigma_i(A) = \lambda_i^{1/2}(AA^*)$ and for a positive definite matrix
 $A$ we have $\sigma_i(A)=\lambda_i(A)$. The trace $\tr(A)$
of a matrix $A$ is given by $\tr(A) = \sum_{i=1}^p \lambda_i(A)$.

We will often be using the inequality
$$   \sum_{i=1}^p \sigma_i^q(AB) \leq \sum_{i=1}^p \sigma_i^q(A)
                         \sigma_i^q(B) ,   \tag 5.1 $$
(see, e.g., Thm.~3.3.14 on p.~176 in \cite{4}),
and certainly we have for a positive definite matrix $$   \| A\|_2 = \max_{1 \leq i \leq p} \lambda_i(A) \leq \tr(A). \tag 5.2 $$

\demo{Proof of $\Rightarrow$}
The proof consists of a few steps, each of which gives an estimate
of an integral.

First we start with
$$   \align  
     \int_0^{2\pi} & \|\prphi - I\|_2 \, d\theta  \\
    &\leq \int_0^{2\pi} \|\left( \prphi \right)^{1/2} - I \|_2 \\
    &\qquad \times \| \left( \prphi \right)^{1/2} + I \|_2 \, d\theta \\
    &\leq 
  \left(  \int_0^{2\pi} \|\left( \prphi \right)^{1/2} - I \|_2^2 \,d\theta 
    \right. \\ 
    &\qquad \times \left.
      \int_0^{2\pi} \|\left( \prphi \right)^{1/2} + I \|_2^2 \, d\theta 
   \right)^{1/2},
 \endalign    $$
where we have used Cauchy-Schwarz for the last inequality. For the second
integral on the right we have by \thetag{5.2}
$$ \aligned
 \int_0^{2\pi} & \|\left( \prphi \right)^{1/2} + I \|_2^2 \, d\theta \\
     &\leq \int_0^{2\pi} \tr \left[ \left( \prphi \right)^{1/2} + I       \right]^2 \, d\theta \\
     &= \int_0^{2\pi} \tr \left[ \prphi  \right. \\
      &\qquad  \left. +\ 2 \left(\prphi \right)^{1/2}
          + I \right]\, d\theta \\
     &= 2\pi p + 2 \int_0^{2\pi} \tr\left( \prphi \right)^{1/2} \, d\theta
        + 2\pi p, 
 \endaligned \tag 5.3  $$
where we used the finite orthogonality \thetag{2.8} for the first term.
To simply the notation, we let $A= \prphi$, then we have for the second term 
$$ \align
 \int_0^{2\pi} \tr\, A^{1/2} \, d\theta
   &= \int_0^{2\pi} \sum_{i=1}^p 
      \lambda_i^{1/2}(A) \, d\theta, \\
   &\leq p^{1/2} \int_0^{2\pi} \left[ \sum_{i=1}^p \lambda_i(A) \right]^{1/2} \, d\theta \\
   &\leq (2\pi p)^{1/2} \left( \int_0^{2\pi} \sum_{i=1}^p \lambda_i(A) \, d\theta \right)^{1/2}
 \endalign $$
where we have used Cauchy-Schwarz for integrals in the last inequality
and Cauchy-Schwarz for sums in the inequality before. We therefore
find, using the finite
orthogonality \thetag{2.8},
$$  \int_0^{2\pi} \tr \left( \prphi \right)^{1/2} \, d\theta
   \leq 2\pi p.  \tag 5.4 $$
Insert this in \thetag{5.3} to find
$$  \multline 
 \int_0^{2\pi} \|\prphi - I\|_2 \, d\theta \\
    \leq (8\pi p)^{1/2} \left( \int_0^{2\pi} \|\left( \prphi \right)^{1/2} - I \|_2^2 \,d\theta \right)^{1/2} . 
 \endmultline   \tag 5.5 $$
Our goal will now be to show that the integral on the right tends to
0 as $n \to \infty$, uniformly for all $\ell \geq 1$. Let again
$A=\prphi$, then we have for this integral
$$ \align
  \int_0^{2\pi} \| A^{1/2} - I \|_2^2 \,d\theta  
 &= \int_0^{2\pi} \| A - 2 A^{1/2} + I \|_2\, d\theta \\
 &\leq \int  \int_0^{2\pi} \tr(A-2A^{1/2}+I) \, d\theta \\
 &\leq 2\pi p - 2 \int_0^{2\pi} \tr\, A^{1/2} \, d\theta + 2\pi p. 
  \endalign  $$
Since we already have the inequality \thetag{5.5}, our goal is to show
$$  \liminf_{n \to \infty} \inf_{\ell \geq 1}
   \int_0^{2\pi} \tr \left( \prphi \right)^{1/2} \, d\theta \geq 2\pi p.
\tag 5.6 $$

Our second step consists in estimating the integral in \thetag{5.6}.
For this we consider the integral
$$   \int_0^{2\pi} \tr [\varphi_n(z) \rho'(\theta) \varphi_n^*(z)]^{1/2}
    \, d\theta = \int_0^{2\pi} \sum_{i=1}^p \sigma_i^{1/2}(\varphi_n(z)
   \rho'(\theta)\varphi_n^*(z)) \, d\theta.  $$
If we write
$$   \varphi_n(z) \rho'(\theta)\varphi_n^*(z) =
   \varphi_n(z) \varphi_{n+\ell}^{-1}(z) 
   \varphi_{n+\ell}(z)\rho'(\theta)\varphi_n^*(z), $$
and then use \thetag{5.1} with $q=1/2$, $A=\varphi_n(z)\varphi_{n+\ell}^{-1}(z)$,
and $B=\varphi_{n+\ell}(z)\rho'(\theta)\varphi_n^*(z)$, then we have
$$  \align
\int_0^{2\pi} & \tr [\varphi_n(z) \rho'(\theta) \varphi_n^*(z)]^{1/2}
    \, d\theta \\
   &\leq \int_0^{2\pi} \sum_{i=1}^p 
   \sigma_i^{1/2}(\varphi_n(z)\varphi_{n+\ell}^{-1}(z)) \
   \sigma_i^{1/2}(\varphi_{n+\ell}(z)\rho'(\theta)\varphi_n^*(z)) 
    \, d\theta \\
   &\leq \int_0^{2\pi} \left( \sum_{i=1}^p
       \sigma_i(\varphi_n(z)\varphi_{n+\ell}^{-1}(z)) 
    \sum_{i=1}^p \sigma_i(\varphi_{n+\ell}(z)\rho'(\theta) \varphi_n^*(z)) \right)^{1/2}
    \, d\theta \\
   & \leq \left( \int_0^{2\pi} \tr \left( \prphi \right)^{1/2} \, d\theta 
    \right. \\
   &\qquad \times \left.
   \int_0^{2\pi}  \sum_{i=1}^p \sigma_i(\varphi_{n+\ell}(z)\rho'(\theta) \varphi_n^*(z)) \, d\theta \right)^{1/2},
\endalign  $$
where we have used Cauchy-Schwarz for integrals in the last inequality
and Cauchy-Schwarz for sums in the inequality before it.
The last integral can be estimated by using \thetag{5.1} with $q=1$,
$A= \varphi_{n+\ell}(z)\rho'(\theta)^{1/2}$ and $B=\rho'(\theta)^{1/2}
\varphi_n^*(z)$:
$$     \align
 \int_0^{2\pi} &  \sum_{i=1}^p 
\sigma_i(\varphi_{n+\ell}(z)\rho'(\theta) \varphi_n^*(z)) \, d\theta \\
  &\leq \int_0^{2\pi} \sum_{i=1}^p \lambda_i^{1/2}(\varphi_{n+\ell}(z)
   \rho'(\theta)\varphi_{n+\ell}^*(z)) \lambda_i^{1/2}
   (\varphi_n(z) \rho'(\theta) \varphi_n^*(z)) \, d\theta \\
  &\leq  \int_0^{2\pi} \left[ \tr(\varphi_{n+\ell}(z)
\rho'(\theta)\varphi_{n+\ell}^*(z))\, \tr(\varphi_n(z)\rho'(\theta)
   \varphi_n^*(z))\right]^{1/2} \, d\theta \\
  &\leq \left( \int_0^{2\pi} \tr(\varphi_{n+\ell}(z)
\rho'(\theta)\varphi_{n+\ell}^*(z))\, d\theta \right.  \\
  &\qquad \times \left. \int_0^{2\pi} \tr(\varphi_{n}(z)
\rho'(\theta)\varphi_{n}^*(z))\, d\theta \right)^{1/2}.  
 \endalign  $$
If we use the finite orthogonality \thetag{2.8}, then, this gives
$$  \int_0^{2\pi}  \sum_{i=1}^p 
\sigma_i(\varphi_{n+\ell}(z)\rho'(\theta) \varphi_n^*(z)) \, d\theta
  \leq 2\pi p . $$
Hence we get
$$ \multline
  2\pi p \int_0^{2\pi} \tr \left( \prphi \right)^{1/2} \, d\theta  \\
   \geq 
   \left(  \int_0^{2\pi} \tr [\varphi_n(z) \rho'(\theta) \varphi_n^*(z)]^{1/2}
    \, d\theta \right)^2.  
 \endmultline   $$
Observe that the right hand side does not depend on $\ell$, hence
$$ \multline
  \inf_{\ell\geq 1} 2\pi p \int_0^{2\pi} \tr \left( \prphi \right)^{1/2} \, d\theta  \\
   \geq 
   \left(  \int_0^{2\pi} \tr [\varphi_n(z) \rho'(\theta) \varphi_n^*(z)]^{1/2}
    \, d\theta \right)^2.  
 \endmultline \tag 5.7  $$

The third step is to estimate the integral on the right. For this
we consider the integral
$$  \int_0^{2\pi} \tr (f(\theta) \rho'(\theta) f^*(\theta))^{1/4}\,
   d\theta= \int_0^{2\pi} \sum_{i=1}^p \sigma_i^{1/2}
  (f(\theta)\rho'(\theta)^{1/2})\, d\theta, $$
where $f$ is a suitable function. Use \thetag{5.1} with $q=1/2$,
$A=f(\theta)\varphi_(z)^{-1}$ and $B=\varphi_n(z)\rho'(\theta)^{1/2}$,
then
$$ \align
  \int_0^{2\pi} \tr (f \rho' f^*)^{1/4}\,
   d\theta 
   &\leq \int_0^{2\pi} \sum_{i=1}^p \sigma^{1/2}(f\varphi_n^{-1})
        \sigma^{1/2}(\varphi_n\rho'(\theta)^{1/2})\, d\theta \\
   &\leq \int_0^{2\pi} \left( \sum_{i=1}^p \sigma(f\varphi_n^{-1})
    \sum_{i=1}^p \sigma(\varphi_n\rho'(\theta)^{1/2})
      \right)^{1/2} \, d\theta \\
   &\leq \left( \int_0^{2\pi} \sum_{i=1}^p \sigma(f\varphi_n^{-1})
     \, d\theta \int_0^{2\pi} \sum_{i=1}^p
        \sigma(\varphi_n\rho'(\theta)^{1/2})\, d\theta \right)^{1/2} \\
   &\leq \left( p \int_0^{2\pi} \sum_{i=1}^p \sigma^2(f\varphi_n^{-1})\,
    d\theta \right)^{1/4} \left( \int_0^{2\pi} \tr(\varphi_n\rho'\varphi_n^*)^{1/2}\, d\theta) \right)^{1/2} \\
  &= \left( 2\pi p \frac{1}{2\pi} \int_0^{2\pi} \tr(f\varphi_n^{-1}
     (\varphi_n^*){-1}f^*)\, d\theta \right)^{1/4} \\
  &\qquad \times 
    \left( \int_0^{2\pi} \tr(\varphi_n\rho'\varphi_n^*)^{1/2}\, d\theta) \right)^{1/2},
\endalign $$
where we use the Cauchy-Schwarz inequality for sums twice
and for integrals once.
if we take the limit as $n \to \infty$, then 
$$ \multline \int_0^{2\pi} \tr (f \rho' f^*)^{1/4}\, d\theta \\
  \leq \left( 2\pi p \tr \int_0^{2\pi} f\, d\rho\, f^* \right)^{1/4} \liminf_{n \to \infty} \left( \int_0^{2\pi} \tr
   (\varphi_n \rho' \varphi_n^*)^{1/2} \, d\theta \right)^{1/2}. \endmultline \tag 5.8 $$
The first integral on the right is
$$ \int_0^{2\pi} f\, d\rho\, f^* = \int_0^{2\pi} f \rho' f^* \, d\theta
  + \int_0^{2\pi} f\, d\rho_s\, f^*, $$
where $\rho_s$ is the singular part of the measure. If $f$ is a matrix function such that $f f^* \leq c I$, 
with $c > 0$, then
$$ \tr \int_0^{2\pi} f G_n\, d\rho_s\, G_n^* f^* \leq
    c \tr \int_0^{2\pi} G_n\, d\rho_s\, G_n^* , $$
and this converges to 0. Therefore we replace $f$ in \thetag{5.8}
by $f G_m$ and let $m \to \infty$ to find
$$ \multline \int_0^{2\pi} \tr (f \rho' f^*)^{1/4}\, d\theta \\
  \leq \left( 2\pi p \tr \int_0^{2\pi} f \rho' f^* \right)^{1/4} \liminf_{n \to \infty} \left( \int_0^{2\pi} \tr
   (\varphi_n \rho' \varphi_n^*)^{1/2} \, d\theta \right)^{1/2}. \endmultline \tag 5.9 $$ 
Now consider the Schur decomposition of $\rho'$
$$   \rho'(\theta) = U \pmatrix
              \lambda_1   & & & \\
		    & \lambda_2 & & \\
	           & & \ddots &  \\
			& & & \lambda_p
		 \endpmatrix  U^*, $$
where $U$ is a unitary matrix function. The hypothesis that
$\det \rho' > 0$ almost everywhere implies that each eigenvalue
$\lambda_i > 0$ almost everywhere. By Lusin's theorem there is
a sequence of nonnegative continuous functions $f_{n,i}$ with
$0 \leq f_{n,i}(\theta \leq 1/\epsilon$, such that
$$  \lim_{n \to \infty}	f_{n,i}(\theta) = \frac{1}{\lambda_i(\theta) + \epsilon}, \qquad \text{in measure}. $$
If we take
$$  F_m = \pmatrix
		f_{m,1}^{1/2} & & & \\
		  & f_{m,2}^{1/2} & &  \\
	        &  & \ddots & \\
		  & & & f_{m,p}^{1/2} \\
		\endpmatrix   U^*, $$
then $F_m F_m^* \leq e^{-1}I$ and 
$$\lim_{m \to \infty} F_m(\theta) \rho'(\theta) F_m^*(\theta) =
   \pmatrix 
      \frac{\lambda_1}{\lambda_1+\epsilon} & & & \\
	 & & \ddots & 
       & &	\frac{\lambda_p}{\lambda_p+\epsilon}
     \endpmatrix := \Lambda_\epsilon \qquad \text{in measure}, $$
hence replacing $f$ by $F_m$ in \thetag{5.9} and using dominated convergence (for convergence in measure) gives
$$ \int_0^{2\pi} \tr \Lambda_\epsilon^{1/4} \, d\theta 
\leq \left( 2\pi p \int_0^{2\pi} \tr\lambda_\epsilon \, d\theta \right)^{1/4} \liminf_{n \to \infty} \left( \int_0^{2\pi} \tr
   (\varphi_n \rho' \varphi_n^*)^{1/2} \, d\theta \right)^{1/2}. \tag 5.10 $$
Finally we let $\epsilon \to 0$, then
$$  \Lambda_\epsilon \to I \qquad \text{almost everywhere}, $$
and by using dominated convergence (for almost everywhere convergence)
we have from \thetag{5.10}
$$ 	2\pi p \leq (2\pi p)^{1/2} 
\liminf_{n \to \infty} \left( \int_0^{2\pi} \tr
   (\varphi_n \rho' \varphi_n^*)^{1/2} \, d\theta \right)^{1/2} $$
which gives the inequality 
$$  \liminf_{n \to \infty} \int_0^{2\pi} \tr
   (\varphi_n \rho' \varphi_n^*)^{1/2} \, d\theta \geq 2\pi p. $$
Combine this with \thetag{5.7}, then the required inequality 
\thetag{5.6} indeed follows. \qed
\enddemo

\enddocument